\documentclass{article}

\usepackage{amsmath,amssymb,amscd}

\newtheorem{Thm}{Theorem}
\newtheorem{Lem}{Lemma}

\newtheorem{Claim}{Claim}

\begin{document}
\title{A comparison between Hofer's metric and ${C^1}$-topology}
\author{Yoshihiro Sugimoto}
\date{}
\maketitle
\begin{abstract}
Hofer's metric is a bi-invariant metric on Hamiltonian diffeomorphism groups. Our main result shows that the topology induced from Hofer's metric is weaker than ${C^1}$-topology if the symplectic manifold is closed.
\end{abstract}

\section{Introduction}
Let ${(M,\omega)}$ be a closed symplectic manifold. In this paper, we define the Hamiltonian vector field of ${H\in C^{\infty}(S^1\times M)}$ by
\begin{equation*}
\omega(X_H,\cdot)=-dH .
\end{equation*}
We denote the time $1$ flow of the vector field ${X_H}$ by ${\phi_H^t}$ and the time $1$ flow by ${\phi_H}$. ${\phi_H}$ is called Hamiltonian diffeomorphism generated by $H$ and we define the Hamiltonian diffeomorphism group ${\textrm{Ham}(M,\omega)}$ by
\begin{equation*}
\textrm{Ham}(M,\omega)=\{\phi_H \ | \ H\in C^{\infty}(S^1\times M)\}.
\end{equation*}
Hofer's norm of a Hamiltonian function ${H\in C^{\infty}(S^1\times M)}$ is defined by
\begin{equation*}
||H||=\int^1_0\max H_t-\min H_tdt .
\end{equation*}
Then, Hofer' norm of ${\phi\in\textrm{Ham}(M,\omega)}$ is defined by
\begin{equation*}
||\phi||=\inf_{\phi_H=\phi}||H||
\end{equation*}
and Hofer's metric on ${\textrm{Ham}(M,\omega)}$ is defined by
\begin{equation*}
\rho(\phi,\psi)=||\phi \psi^{-1}||.
\end{equation*}

On the Hamiltonian diffeomorphism group ${\textrm{Ham}^c(M,\omega)}$, we can considered ${C^1}$-topology and the topology induced from Hofer's metric. In this paper, we compare these two topologies. We denote the families of all open subsets of these two topologies by ${\mathcal{O}_{C^1}}$ and ${\mathcal{O}_{\mathrm{Hofer}}}$. We prove the following theorem.
\begin{Thm}
Let ${(M,\omega)}$ be a closed symplectic manifold. Then, ${C^1}$-topology is stronger than the topology induced by Hofer's metric. In other words, ${\mathcal{O}_{C^1}\subset \mathcal{O}_{\mathrm{Hofer}}}$ holds.
\end{Thm}

\section{Proof of the main result}
It suffices to prove the following claim.
\begin{Claim}
Let ${V_{\mathrm{Hofer}}\in \mathcal{O}_{\mathrm{Hofer}}}$ be an open subset which contains the identity. Then, there is another open neighborhood ${U_{C^1}\subset \mathcal{O}_{C^1}}$ of the identity so that ${U_{C^1}\subset V_{\mathrm{Hofer}}}$ holds.
\end{Claim}

We use the following correspondence between symplectomorphisms which are ${C^1}$-close to the identity and closed ${1}$-forms on $M$. Let ${\sigma \in \Omega^2(M\times M)}$ be a symplectic form defined by ${\sigma=(-\omega)\oplus \omega}$ and let ${\lambda \in \Omega^1(T^*M)}$ be the canonical Liouville $1$-form on the cotangent bundle. We can choose a neighborhood ${\mathcal{N}(M_0)\subset T^*M}$ of the zero section ${M_0}$ and a neighborhood ${\mathcal{N}(\Delta)}\subset M\times M$ of the diagonal ${\Delta \subset M\times M}$ and a symplectomorphism ${\Psi}$ as follows.

\begin{gather*}
\Psi: \mathcal{N}(\Delta)\longrightarrow \mathcal{N}(M_0) \\
\Psi(\Delta)=M_0  \\
\pi(\Psi(q,q))=q
\end{gather*}

Let ${\psi\in \mathrm{Ham}(M,\omega)}$ be sufficiently ${C^1}$-close to the identity. Then, the corresponding closed $1$-form ${\sigma_{\psi}}$ is defined by
\begin{equation*}
\sigma_{\psi}=\Psi(\mathrm{graph}(\psi))
\end{equation*}

We fix a Hamiltonian function ${H\in C^{\infty}(S^1\times M)}$ so that ${\phi_H=\psi}$ holds and a path ${\psi_t}$ of symplectomorphisms as follows.
\begin{equation*}
(q,\psi_t(q))=\Psi^{-1}(t\sigma_{\psi})
\end{equation*} 
In other words, ${\sigma_{\psi_t}=t\sigma_{\psi}}$ holds. 
\begin{Lem}
This path satisfies ${\psi_t\in \mathrm{Ham}(M,\omega)}$.
\end{Lem}
For this lemma, it suffices to prove that ${\sigma_{\psi}}$ is an exact $1$-form. We choose a loop ${l_t}$ of symplectomorphisms as follows.
\begin{equation*}
l_t=\begin{cases} \psi_{2t} & 0\le t\le \frac{1}{2} \\ \phi_H^{2-2t}  & \frac{1}{2}\le t \le 1  \end{cases}
\end{equation*}
Then, the flux homomorphism of this path can be calculated as follows (see \cite{MS} Lemma 10.15).
\begin{gather*}
\mathrm{Flux}(\{l_t\})=[\sigma_{\psi}]\in H^1(M:\mathbb{R})
\end{gather*}
The flux group ${\Gamma\subset H^1(M:\mathbb{R})}$ is known to be discrete (\cite{O}). So, by making ${\psi}$ sufficiently ${C^1}$-close to the identity, ${[\sigma_{\psi}]}$ becomes ${0\in H^1(M:\mathbb{R})}$. This implies that ${t\sigma_{\psi}}$ are exact and ${\psi_t}$ are Hamiltonian diffeomorphisms. So ${\psi_t}$ is a path of Hamiltonian diffeomorphisms between the identity and ${\psi}$ and we proved Lemma 1.  \begin{flushright}    $\Box$ \end{flushright}

Let $K\in C^{\infty}(S^1\times M)$ be a Hamiltonian function so that ${\phi_K^t=\psi_t}$ holds. By making the path $\psi_t$ ${C^1}$-small, we can make the Hofer norm ${||K||}$ becomes arbitrary small. This implies that by making the neighborhood ${U_{C^1}}$ small,  we can make the Hofer diameter of ${U_{C^1}}$ arbitrary small. So we can make ${U_{C^1}\subset V_{\textrm{Hofer}}}$ and we proved Claim 1.
\begin{flushright}    $\Box$ \end{flushright}

\end{document}